\documentclass[11pt,a4paper]{article}
\usepackage[utf8]{inputenc}
\usepackage[T1]{fontenc}
\usepackage{amssymb}
\usepackage{amsthm}

\usepackage{color}
\definecolor{DarkOlive}{rgb}{0.1047,0.2412,0.0064}
\definecolor{FireBrick}{rgb}{0.5812,0.0074,0.0083}
\definecolor{RoyalBlue}{rgb}{0.0236,0.0894,0.6179}
\definecolor{RoyalGreen}{rgb}{0.0236,0.6179,0.0894}
\definecolor{RoyalRed}{rgb}{0.6179,0.0236,0.0894}
\definecolor{LightBlue}{rgb}{0.8544,0.9511,1.0000}
\definecolor{Black}{rgb}{0.0,0.0,0.0}
\definecolor{MidnightBlue}{rgb}{0.0035,0.0020,0.1363}
\definecolor{FireBrick3}{rgb}{0.5812,0.0074,0.0083}
\definecolor{FireBrick4}{rgb}{0.2156,0.0023,0.0035}
\definecolor{Blue2}{rgb}{0.0000,0.0000,0.8644}
\definecolor{Navy}{rgb}{0.0000,0.0000,0.1927}
\definecolor{MediumBlue}{rgb}{0.0000,0.0000,0.6179}
\usepackage[
        bookmarks=false,pdftitle={Green Functions in Small Characteristic},
        colorlinks=true,backref=false,breaklinks=true,linkcolor=MediumBlue,
        citecolor=FireBrick3,filecolor=RoyalRed,
        urlcolor=Blue2]{hyperref}

\theoremstyle{plain}
\newtheorem{Thm}{Theorem}[section]

\newtheorem{Pro}[Thm]{Proposition}

\theoremstyle{definition}

\newtheorem{Rem}[Thm]{Remark}

\newenvironment{Prf}{\noindent\textbf{Proof.}}{\hfill $\Box$ \medskip}

\newcommand{\hgt}{ \textrm{ht} }
\newcommand{\F}{\mathbb{F}}
\newcommand{\Fq}{\mathbb{F}_q}

\newcommand{\Z}{\mathbb{Z}}

\newcommand{\C}{\mathbb{C}}
\newcommand{\Q}{\mathbb{Q}}

\newcommand{\bG}{\mathbf{G}}
\newcommand{\bT}{\mathbf{T}}
\newcommand{\bB}{\mathbf{B}}
\newcommand{\bU}{\mathbf{U}}
\newcommand{\bL}{\mathbf{L}}
\newcommand{\bP}{\mathbf{P}}
\newcommand{\cL}{\mathcal{L}}

\newcommand{\GF}[1]{\F_{#1}}
\newcommand{\GFq}{\GF{q}}

\newcommand{\barGF}[1]{\bar{\F}_{#1}}

\begin{document}
\title{Green Functions in Small Characteristic}
\author{Frank Lübeck\thanks{This is a contribution to 
Project-ID 286237555 -- TRR 195 -- by the
German Research Foundation.}}
\maketitle
\begin{abstract}
The values of the ordinary Green functions are known for almost all groups
of Lie type, a long term achievement by various authors.

In this note we solve the last open cases, which are for exceptional groups
$E_8(q)$ where $q$ is a power of $2$, $3$ or $5$.
\end{abstract}

\section{Introduction}\label{secIntro}

Let $p$ be a prime, $q$ some power of $p$ and $\GFq$ the field with $q$
elements. Let $\bG$ be a connected reductive algebraic group over an 
algebraic closure ${\barGF{p}}$ with a Frobenius endomorphism $F$ 
defining an $\GFq$-rational structure. 

We are interested in class functions of the finite group of fixed points
$G(q) := \bG^F$. Let $\bT \subset \bG$ be an $F$-invariant maximal torus. 
Deligne and Lusztig~\cite{DL76}
defined for each irreducible character $\theta$ of the
abelian group $\bT^F$ a generalized character $R_{\bT,\theta}^\bG$ of
$G(q)$, using certain modules constructed by $\ell$-adic cohomology.
They show that the values of $R_{\bT,\theta}^\bG$ on 
unipotent elements only depend on the
$G(q)$-conjugacy class of $\bT$ and not on the
character $\theta$, we write  $Q_\bT^\bG$ for the restriction of 
$R_{\bT,\theta}^\bG$ to unipotent elements, this is called the \emph{Green
function} of $\bT$. 

In~\cite{CSII} Lusztig defined another set of class functions on the
unipotent elements using character sheaves, these are called
\emph{generalized Green functions}. A subset of them is also associated to 
$G(q)$-classes of $F$-invariant maximal tori, and Lusztig~\cite{Lu90} (for
large $q$) and Shoji~\cite[2.2]{Sho95} showed
later that the two types of Green functions coincide.

While from the $\ell$-adic cohomology approach it is not so clear how to
compute the Green functions explicitly, the definition via character sheaves 
leads to an algorithm to compute the (generalized) Green functions,
see~\cite{CSV}.

More  precisely, this  algorithm determines  the Green  functions as  linear
combinations of certain functions which are supported on elements $C^F$ for a
single unipotent class $C$ in $\bG$. The values of these functions are clear
up to a normalization by some scalar in $\mathbb{C}$ of absolute value $1$.

Finding these scalars is a non-trivial task. Computing the (generalized)
Green functions for general $\bG$ can be reduced to the case of simple
simply-connected  groups. All cases of groups of classical types were
systematically considered by Shoji~\cite{Sh06,Sh07,Sh22} (also generalized
Green functions), the cases 
of small rank exceptional groups can be read off from the known character
tables, Green functions of groups of type $E_6$, $E_7$ and $E_8$ in 
good characteristic 
were considered by Beynon and Spaltenstein~\cite{BS84}, and various
exceptional groups in small characteristic by Malle~\cite{M93},
Porsch~\cite{P94} and more recently Geck~\cite{Ge20}. In this paper
we describe a method that enabled us to also handle the only case which was
left open so far, that is the groups of type $E_8$ in bad characteristic. 

So, with the results of this paper the Green functions are known in all
cases. It only remains to consider the other generalized Green functions
in groups of exceptional types. This is done in~\cite{LSh24}, using
the results of this article.

Our method is a variant of the idea in~\cite{Ge20} where the permutation
character of the Borel subgroup in $G(q)$ and computations in a matrix
representation were used. Here, we use more general parabolic subgroups and
compute with the Steinberg presentation of the considered groups.

Explicit tables of Green functions for exceptional groups are available
from~\cite{Lu24}.

\section{Notations}\label{secNot}

For more details on the following basic setup we refer to the introductory
sections of the text books~\cite{Ca85, DM20}.

Let $\bG$ be a connected reductive group over an algebraic closure
$\mathbf{k}$ of a
finite field $\Fq$ with $q$ elements in characteristic $p$. 
We assume that $\bG$ is defined over $\Fq$ and call $F$ the corresponding
Frobenius endomorphism of $\bG$. We write $G(q) := \bG^F$ for the finite group
of $F$-fixed points. 

Let  $\bT \subseteq \bG$ be  an  $F$-stable maximal  torus of  $\bG$ that  is
contained in  an $F$-stable  Borel subgroup $\bB  \subseteq \bG$.  The group
$\bG$ is  determined up  to isomorphism  by its root  datum with  respect to
$\bT$. The minimal unipotent subgroups  $X_r$ normalized by $\bT$ are called
root subgroups.  There is an isomorphism  $\mathbf{k}^+ \to X_r$, $a  \mapsto
x_r(a)$.
The root $r:\bT \to \mathbf{k}^\times$ is an  element 
of the  (additive) character
group $X(\bT)$  and describes  the conjugation: $t x_r(a) t^{-1}   = x_r(r(t)
a)$. We write $\Phi$ for the finite set of roots, then we have
$\Phi = \Phi^+ \dot{\cup} \Phi^-$ where a root $r$ lies in $\Phi^+$ when
$X_r \subseteq \bB$; these subsets are called positive roots and negative
roots, respectively, and we have $\Phi^- = -\Phi^+$.
The positive roots contain a unique subset $\Delta$
such that every positive root is a unique non-negative linear combination
of the roots in $\Delta$; the elements of $\Delta$ are called simple roots.
The height $\hgt(r)$ of a positive root $r \in \Phi^+$ is the sum of its
coefficients when $r$ is written as linear combination of the roots in
$\Delta$.

The group $\bG$ is generated by $\bT$ and all root subgroups 
$X_r$, $r \in \Phi$.
The Borel subgroup $\bB$ is generated by $\bT$ and the $X_r$ with $r \in
\Phi^+$. The group $\bU = \prod_{r \in \Phi^+} X_r$ is a maximal unipotent
subgroup and the unipotent radical of $\bB$. 
We will describe $\bU$ in more detail in
Section~\ref{secUcomp}.

The group $W = N_\bG(\bT)/\bT$ is called the Weyl group of $\bG$, it is a
Coxeter group with Coxeter generators $S := \{s_r\mid\; r \in \Delta\}$, where 
$s_r$ is the unique non-trivial coset of $\bT$ in $N_\bG(\bT) \cap \langle
X_r, X_{-r}\rangle$. The length $l(w)$ of an element $w \in W$ is the smallest
integer $k$ such that $w = s_1 \cdots s_k$ with $s_i \in S$, the sequence 
$s_1, \ldots, s_k$ is called a reduced word for $w$.
An element $w \in W$ permutes the roots $r \in \Phi$ where $w(r)$ is defined
by $w X_r w^{-1} = X_{w(r)}$. For $r \in \Delta$ we have $s_r(r) = -r$ 
and $s_r(\Phi^+ \setminus \{r\}) \subset \Phi^+$. So, for 
$w \in W$ and $r \in \Delta$ we have $l(w s_r) < l(w)$ (and in that case 
$l(w s_r) = l(w)-1$) if and only if $w(r) \in \Phi^-$.

The  Frobenius  endomorphism $F$  restricts  to  Frobenius endomorphisms  of
$\bB$,  $\bT$  and $\bU$  and  induces  a natural  map  on  $X(\bT)$ and  an
automorphism of $W$ which permutes the  set $S$ of generators. We will write
$F$  also for  the induced  maps,  and $H(q)  =\mathbf{H}^F$ for  $F$-stable
subgroups $\mathbf{H}$.

For $w \in W$ let $\dot{w} \in N_\bG(\bT)$ be a representative. Using the
Lang-Steinberg theorem we can find a $g \in \bG$ with $g F(g^{-1}) =
\dot{w}$. Then $\bT_w := \bT^g$ is also an $F$-stable maximal torus of $\bG$
(well defined up to $G(q)$-conjugacy). Let $M = \{w_1,\ldots, w_k\}$ be a
set of representatives for the $F$-conjugacy classes of $W$ (where $w, w'
\in W$ are $F$-conjugate if there is a $z \in W$ such that $w' = z w
F(z^{-1})$), then $\{ \bT_{w_1}, \ldots, \bT_{w_k}\}$ is a set of
representatives for the $G(q)$-conjugacy classes of $F$-stable maximal tori.

We will consider the Bruhat decomposition of $\bG$ as union of disjoint
double cosets with respect to the Borel subgroup $\bB$:
\[ \bG = \dot{\bigcup_{w \in W}} \bB \dot{w} \bB,\]
where $\dot{w} \in N_\bG(\bT)$ is again a representative of $w$.
We will use a stronger form of this Bruhat decomposition that says, that
an element $g \in \bB \dot{w} \bB$ can be uniquely written as product
$g = u t \dot{w} u'$ with $u \in \bU$, $t \in \bT$, $u' \in \bU_w = 
\prod_{r \in \Phi^+, w(r) \in \Phi^-} X_r$.

When $F$ acts trivially on $\Phi$ and $W$ we can choose the $\dot{w} \in
N_\bG(\bT)^F$ and have the same decomposition for
the corresponding finite groups:
\[ G(q) = \dot{\bigcup_{w \in W}} U(q) T(q) \dot{w} U_w(q).\]

We will need one basic property of groups with a $BN$-pair which have such a
Bruhat decomposition: For $w \in W$ and $s \in S$ we have 
\[ (BsB)(BwB) \left\{ \begin{array}{ll}
= B(sw)B, &\textrm{if } l(sw) > l(w)\\
\subseteq (BwB) \cup (B (sw) B),& \textrm{otherwise.}
\end{array} \right. \]

\section{Computing in the unipotent subgroup}\label{secUcomp}

We want to do certain computations with elements in $\bG$ or $G(q)$ where
the elements are represented in the strong form of the Bruhat decomposition.
For this one first needs to fix the 
isomorphisms $\mathbf{k}^+ \to X_r$, $t \mapsto
x_r(t)$ and the
representatives $\dot w$ of Weyl group elements. 
And for the multiplication one needs to know commutator rules for the
factors appearing in Bruhat expressions of the form $g = u t \dot{w} u'$.
These are provided by the relations in the Steinberg presentation of
the group $\bG$ or $G(q)$, see for example~\cite[Ch.9]{Spr98}.

It turns out that for our present applications we only need detailed 
computations in the unipotent subgroup $\bU$.  This needs the most 
complicated part of the Steinberg presentation and depends on certain
choices (leading to isomorphic groups). We follow the construction of $\bG$
as a Chevalley group by Carter in~\cite[Ch.4-8]{Ca72}. This construction starts with
a semisimple Lie algebra $\cL$ over $\C$ which has the same type of root
system as $\bG$. It uses a Chevalley basis of $\cL$ consisting of elements
$h_r$, $r \in \Delta$, spanning a Cartan subalgebra, and of
root vectors $e_r$, $r \in \Phi$. The $e_r$ are acting as nilpotent linear
maps on $\cL$ and can be used to construct the elements in the root subgroups:
\[ x_r(a) := \exp(a \cdot e_r) \quad \textrm{ for $r\in \Phi$, $a \in \C$}.\]
The Chevalley bases have the property that all its structure constants 
(and so the entries of matrices of $e_r$ acting on $\cL$ with respect to 
the Chevalley
basis) are in~$\Z$. So, for each $p$ (prime or $0$) we can reduce the matrix
entries modulo $p$ and use the definition of $x_r(a)$ for any element
$a$ in a field of characteristic $p$. 

\subsection{Different Chevalley bases}\label{extraspecial}
We use Carter's description in~\cite[4.2]{Ca72}.
Let $\mathcal{B} = \{h_s, e_r\mid \; 
s \in \Delta, r \in \Phi\}$ and  $\mathcal{B'} = \{h'_s, e'_r\mid \;
s \in \Delta, r \in \Phi\}$ be two Chevalley bases. Then there is a unique
automorphism of $\cL$ with $h_s \mapsto h'_s$ for $s \in \Delta$,
$e_r \mapsto \lambda_r e'_r$ for $r \in \Phi$ with $\lambda_r = 1$
for $r \in \Delta$, $\lambda_{-r} = \lambda_r \in \{\pm 1\}$ 
for all $r \in \Phi$.

The possibilities for the signs can be described as follows, relative to an
enumeration $\Delta = \{r_1, \ldots, r_l\}$ of the simple roots. Assign to 
each root $r \in \Phi^+ \setminus \Delta$ a sign $\epsilon_r$
by writing $r = \tilde{r} + r_i$
with $\tilde{r} \in \Phi^+$ and $r_i \in \Delta$ for the largest possible $i$
(Carter calls $(\tilde{r},r_i)$ an extraspecial pair), then 
$[e_{\tilde{r}},e_{r_i}] =
N_{\tilde{r},r_i} e_r \neq 0$ and 
$\epsilon_r$ is the sign of $N_{\tilde{r},r_i}$ (this is
up to sign the same number for any Chevalley basis). 
If we also know the signs $\epsilon'_r$ for $r \in \Phi^+$, defined in the
same way with respect to the basis $\mathcal{B}'$, we can consider all roots
$r \in \Phi^+$ by increasing height and determine if the isomorphism above maps
$e_r$ to $e'_r$ or $-e'_r$:  
We have  for $r \in \Phi^+ \setminus \Delta$, $r = \tilde{r} + r_i$ as above,
\[ \lambda_r = \lambda_{\tilde{r}} \epsilon_r \epsilon'_r.\]

This describes the diagonal transition matrix between the two Chevalley bases
and also determines an isomorphism between the groups generated by the root
subgroups $X_r$ or $X'_r$,  constructed with respect to the different bases.
It maps $x_r(a) \mapsto x'_r(\lambda_r a)$ for all $r \in \Phi$.

For any choice of signs for the extraspecial pairs there exists a
corresponding Chevalley basis.

We mention that Geck~\cite{Ge17} described an (up to a global sign) 
canonical choice of a Chevalley basis that does not depend on certain
choices of signs.

\subsection{Commutator formula}\label{comm}
In~\cite[4.2.2]{Ca72} there is a detailed description of how to compute for a
given root system and given signs for extraspecial pairs all (uniquely
determined) structure constants of the corresponding Chevalley basis of the
Lie algebra $\cL$.

Furthermore, using the formulae in~\cite[4.3.1,5.2.2]{Ca72} we can
compute for any
two different roots $r_1, r_2 \in \Phi^+$ and any $i,j \in \Z_{>0}$ such
that $i r_1 + j r_2 \in \Phi^+$ an integer $C_{ijr_1r_2}$. These enable us
to compute with elements in the group $\bU = \langle X_r\mid\; r \in
\Phi^+\rangle$ using the following commutator formula~\cite[5.2.3]{Ca72}:
\[ x_{r_2}(a_2) x_{r_1}(a_1) = x_{r_1}(a_1) x_{r_2}(a_2)  
\prod_{i,j} x_{i r_1 + j r_2}(C_{ijr_1r_2} (-a_1)^i a_2^j),\]
where the product is over all $i,j \in \Z_{>0}$ with $i r_1 + j r_2 \in
\Phi$ sorted by increasing $i+j$ (factors with the same $i+j$ always
commute).

\begin{Pro}\label{reorder}
Let $r_1, r_2, \ldots, r_N$ be the positive roots of $\bG$ in any fixed
order, we write $x_i(a) := x_{r_i}(a)$ for the corresponding root
elements.
\begin{itemize}
\item[(a)] Any element of $\bU$ can be uniquely written in the form
\[ x_1(a_1) x_2(a_2) \cdots x_N(a_N).\]
\item[(b)] Any product of root elements $x_{t_1}(b_1) \cdots
x_{t_k}(b_k)$ where $r_{t_1}, \ldots, r_{t_k}$ are positive roots can 
be rewritten to 
the form in~(a) (where some factors $x_i(0) = 1$ may be omitted) 
by applying a finite number of the following steps to
any pair of consecutive factors $x_{t_i}(b_i)x_{t_{i+1}}(b_{i+1})$ with
$t_{i+1} \leq t_i$: 
\begin{itemize}
\item If $t_i = t_{i+1}$ simplify to one factor $x_{t_i}(b_i+b_{i+1})$.
\item Otherwise substitute the two factors according to the commutator
formula.
\end{itemize}
In practice it works best to handle pairs with commuting factors first.
\item[(c)] For positive roots $r_i, r_j$ we write $r_i \prec r_j$ 
if and only if $r_j - r_i$ is a non-zero non-negative linear combination 
of positive roots.

Let $r_i$ be a positive root, and
\[ u = x_1(a_1) x_2(a_2) \cdots x_N(a_N) x_i(a).\]
After reordering factors as in~(b) we get a product
\[ u = x_1(b_1) x_2(b_2) \cdots x_N(b_N).\] 
Then we have $b_i = a_i + a$ and $b_j = a_j$ whenever $j \neq i$ and
$r_i \not\prec r_j$.

\end{itemize}
\end{Pro}
\begin{Prf}
To show~(b) we have to show that the process terminates. For this we count
for each positive integer $m$ the number of pairs $(t_i,t_j)$ with $i < j$,
$t_i > t_j$ and $\hgt(r_{t_1}) + \hgt(r_{t_j}) = m$. It is clear that after
the first type of substitution in~(b) none of these counts will be larger
than before. The substitution using the commutator formula will reduce the
count for $m = \hgt(r_{t_i}) + \hgt(r_{t_{i+1}})$ by one and maybe enlarge
the counts for some larger $m$. Since there is an upper
bound for the height of all roots, all counts will be zero after a finite
number of steps, so the factors are sorted.

Statement~(c) follows by applying~(b) to the given $u$ and noticing that the
commutator formula applied to $x_{r}(a)x_{r'}(a')$ only introduces new
factors $x_{r''}(a'')$ with $r \prec r''$ and $r' \prec r''$.

Part~(b) yields a constructive proof of the existence of the form in~(a).
The uniqueness of the factorization is shown in~\cite[5.3.3]{Ca72} in the
case that the ordering of the roots refines the ordering by height.
The general case follows by induction over the lowest height of roots
$r_i$ with non-trivial factor $x_i(b_i)$; as in~(c) those factors will
not change when moved to the left. \mbox{}
\end{Prf}

\section{Green functions by the Lusztig-Shoji algorithm}\label{secGFalg}

Deligne and Lusztig~\cite{DL76} defined for each $F$-stable maximal torus
$\bT$ and each irreducible character $\theta$ of the abelian group $T(q)$ 
a generalized character $R_{\bT,\theta}^\bG$ of $G(q)$. Their values on 
unipotent elements depend only on the $G(q)$-conjugacy class of $\bT$ and
not on $\theta$. The restrictions $Q_\bT^\bG$ to unipotent elements are 
called the (ordinary) Green functions of $G(q)$. 

As in Section~\ref{secNot} let $\{ \bT_{w_1}, \ldots, \bT_{w_k}\}$ be a set
of representatives of the $G(q)$-conjugacy classes of maximal tori, where
$\{w_1,\ldots, w_k\} \subset W$ are representatives of the $F$-conjugacy
classes of $W$. We write $Q_{w_i} := Q_{\bT_{w_i}}^\bG$. 

Let $\chi \in \textrm{Irr}(W)$ be an $F$-stable irreducible character of
$W$. Writing $F$ also for the automorphism of $W$ induced by the Frobenius
endomorphism, we consider the semidirect product $W \rtimes \langle F
\rangle$; the $W$-conjugacy classes in the coset $WF \subseteq W \rtimes
\langle F \rangle$ (the $W$-coset containing $F$) are in bijection with the
$F$-conjugacy classes of $W$. The character $\chi$ can be extended to an
irreducible character $\tilde\chi$ of $W \rtimes \langle F \rangle$.

For any class function $\tilde\chi$ on the $W$-conjugacy classes of the
coset $WF$ we define the linear combination

\[Q_{\tilde\chi}  :=  
\frac{1}{|W|} \sum_{w \in W} \tilde\chi(wF) Q_w =
\sum_{i=1}^k \frac{1}{|C_W(w_i F)|} \tilde\chi(w_i F) Q_{w_i}.\]

Let $C$ be an $F$-stable unipotent class of $\bG$ and $u \in C^F$.
The Lang-Steinberg theorem shows that the $G(q)$-conjugacy classes in 
$C^F$ are parameterized by the $F$-conjugacy classes of the component group
$A(u) = C_\bG(u)/C_\bG^0(u)$ (the centralizer of $u$ in $\bG$ modulo its
connected component), or the $A(u)$-conjugacy classes in the coset
$A(u)F \subseteq A(u) \rtimes \langle F \rangle$ where now $F$ denotes the
automorphism on $A(u)$ induced by the Frobenius endomorphism. We write $u_a$
for an element in the $G(q)$-class of $C^F$ corresponding to $a \in A(u)$.
An $F$-stable irreducible character $\rho \in \textrm{Irr}(A(u))$ can be
extended to an irreducible character $\tilde\rho$ of
$A(u) \rtimes \langle F \rangle$. We consider the following class function
on $G(q)$:
\[ Y_{u,\tilde\rho}(g) := \left\{
\begin{array}{ll}
\tilde\rho(aF),& \textrm{ if $g$ is conjugate to $u_a$}\\
0, & \textrm{ else.}
\end{array}
\right. \]
Using the theory of character sheaves Lusztig~\cite[24.]{CSV} 
described an algorithm to write another set of class functions, also denoted
$Q_{\tilde\chi}$, as linear combinations of
the functions $\zeta_{u,\tilde\rho} Y_{u,\tilde\rho}$, where
$\zeta_{u,\tilde\rho} \in \C$ are scalars of absolute value $1$ which are
not determined by the algorithm. The main input  of the algorithm 
is the Springer correspondence which is an injective map from
$\textrm{Irr}(W)$ to the set of pairs $(u,\rho)$ modulo $\bG$-conjugacy.
Further data which are needed are the dimensions of the unipotent classes of
$\bG$ and the $F$-character table of the Weyl group.

Later, Lusztig showed in~\cite{Lu90} that the two types of class functions
denoted $Q_{\tilde\chi}$ coincide under some conditions on $q$, and Shoji
showed in~\cite{Sho95} that the same holds without restrictions on $q$.

For much more detailed descriptions of this setup  we refer 
to~\cite{CSV,Sh06,Ge20a,Ge20}.

The following proposition reduces the determination of the $\zeta_{u,\rho}$
for untwisted groups to the special case $q = p$, that is the groups $G(p)$ 
defined over the prime field. It is a special case of~\cite[Thm.
3.7]{Ge20a}.
\begin{Pro}\label{proprimefield}
Assume that the action of $F$ on $\Phi$ is trivial.
Let $C$ be an F-stable unipotent class with representative $u \in C^F = C
\cap G(q)$. Assume that $F$ acts trivially on $A(u)$ (so that we have
$\tilde\rho = \rho$ in the discussion above). Let $\rho \in \textrm{Irr}(A(u))$
appear in the Springer correspondence and $\zeta_{u,\rho}$ be the associated
scalar. Then we have for any positive integer $m$ that $u$ is also stable 
under $F^m$. Write $\zeta^{(m)}_{u,\rho}$ for the scalar associated to
$(u,\rho)$ when $u$ is considered as element of $\bG^{F^m} = G(q^m)$. Then
we have $\zeta^{(m)}_{u,\rho} = (\zeta_{u,\rho})^m$.
\end{Pro}

What remains to find the Green functions $Q_w$ (or equivalently the
$Q_{\tilde\chi}$ for $F$-stable irreducible $\chi$) 
as class functions is to choose 
representatives $u$ of unipotent classes and the characters $\tilde\rho$ and
then to determine the scalars $\zeta_{u,\tilde\rho}$. As mentioned in
Section~\ref{secIntro} this task has been done in almost all cases.
We describe a method which enabled us to handle the remaining cases of
groups $E_8(q)$ in bad characteristic $2$, $3$ and $5$.

For computations we used our own implementation of the Springer
correspondence and the Lusztig-Shoji algorithm, but all that is needed is
also available in Michel's version of CHEVIE~\cite{JCHEVIE}.

The following additional information is also useful, see~\cite[24.]{CSV}.
\begin{Rem}\label{remLambda}
The Lusztig-Shoji algorithm also returns the matrix $\Lambda$ whose rows and
columns are labeled by the pairs $(u,\tilde\rho)$ as above and the entry
in position $(u,\tilde\rho)$, $(u',\tilde\rho')$ is
\[ \sum_{v \in G(q)\textrm{ \scriptsize unipotent} }
\zeta_{u,\tilde\rho} Y_{u,\tilde\rho}(v) 
\overline{\zeta_{u',\tilde\rho'} Y_{u',\tilde\rho'}(v)}.\]
This matrix is block diagonal (one block for each unipotent $F$-stable 
class in $\bG$), it is
symmetric and has values in the rational numbers.

From this we can always determine the sizes of the $G(q)$-classes in $C^F$ for
each class: If for a class $C$ there are $k$ classes in $C^F$ and we write $Y_C$
for the $k \times k$ matrix whose rows contain the values of 
the $\zeta_{u,\tilde\rho} Y_{u,\tilde\rho}$, and if 
$\Lambda_C$ is the corresponding diagonal block 
of $\Lambda$, then $Y_C^{-1} \Lambda_C \overline{Y_C}^{-t}$ is 
a diagonal matrix where the diagonal entries are the class lengths.

These properties sometimes yield restrictions on the values of the
scalars $\zeta_{u,\tilde\rho}$.
\end{Rem}

\section{Permutation characters of parabolic subgroups}\label{secP}

For background about parabolic subgroups we refer to~\cite[Ch.3]{DM20}.

We want to consider standard parabolic subgroups of $G(q)$. These are
parameterized by $F$-stable subsets $J \subseteq \Delta$ of the simple
roots. The subgroup  $W_J \leq W$ generated by the set 
$S_J  := \{s_r\mid\; r \in J\}$ is also
a Coxeter group with $S_J$ as set of Coxeter generators.
The set $\bP_J = \bigcup_{w \in W_J} \bB \dot{w} \bB$ is a subgroup of $\bG$
and is called a standard parabolic subgroup of $\bG$, and similarly for the
finite groups, that is $P_J(q) = \bigcup_{w \in W_J} B(q) \dot{w} U_w(q) 
\leq G(q)$. 

The subset $\Phi_J := W_J(J) \leq \Phi$ is also a root system and the
parabolic subgroup has a Levi decomposition  $\bP_J = \bL_J \bU_J$, where 
$\bU_J$ is
the unipotent radical, generated by the $X_r$ with $r \in \Phi^+ \setminus 
\Phi_J$, and $\bL_J$ is a Levi complement, it is generated by $\bT$ and the
root subgroups $X_r$ with $r \in \Phi_J$.

\subsection{Permutation characters by Deligne-Lusztig
characters}\label{ssecPDL}

We show how to find the values of permutation characters of parabolic
subgroups as linear combinations of Deligne-Lusztig characters
$R_{\bT_w,1}^\bG$.

Recall that we have a set of representatives $\{w_1,\ldots,w_k\}$ of the
$F$-conjugacy classes of $W$ and that the $Q_{w_i}$ are the corresponding Green
functions.

First we mention that the trivial character $1_{G(q)} \in \textrm{Irr}(G(q))$
is a linear combination of
Deligne-Lusztig characters~\cite[10.2.5]{DM20}:
\[ 1_{G(q)} =  \sum_{i=1}^k \frac{1}{|C_W(w_i F)|}  R_{\bT_{w_i},1}. \]
So, its restriction to unipotent elements is $Q_{\bar\chi}$ for the trivial
character $\bar\chi$ of the coset $WF$ (see the definition of
$Q_{\tilde\chi}$ in Section~4).

The trivial character on a parabolic subgroup $P_J(q)$ is the inflation of the
trivial character of $L_J(q)$ to $P_J(q)$ via the canonical map 
$P_J(q) \to L_J(q)$.
Therefore, the permutation character of $P_J(q) \leq G(q)$ is the
Harish-Chandra induction of the trivial character on $L_J(q)$ which is a
special case of Lusztig induction $R_{\bL_J}^\bG$, see~\cite[5.]{DM20}.
Let $\{v_1, \ldots, v_m\}$ be representatives of the $F$-conjugacy classes in
$W_J$.
Then we have $1_{L_J(q)} =  \sum_{j=1}^m (1/|C_{W_J}(v_j F)|)
R_{\bT_{v_j},1}^{\bL_J}$.
Using the transitivity of Lusztig induction~\cite[9.1.8]{DM20} we get
\[ {1_{P_J(q)}}^{G(q)} = 
R_{\bL_J}^\bG (\sum_{j=1}^m \frac{1}{|C_{W_J}(v_j F)|} 
R_{\bT_{v_j},1}^{\bL_J}) = 
\sum_{j=1}^m (\frac{1}{|C_{W_J}(v_j F)|} R_{\bT_{v_j},1}^{\bG}).\]
This shows that the restriction to unipotent elements is $Q_{\bar\chi}$
where $\bar\chi$ is the permutation character of $WF$ on $W_JF \subseteq WF$
(defined as $\bar\chi(wF) =1/|W| \cdot \sum_{x \in W,\; (wF)^x \in W_JF} 1$, 
which specializes to ${1_{W_J}}^W$ in case $F = \textrm{id}_W$).

Given all $Q_{w_i}$, where the values are written with the so far unknown
scalars $\zeta_{u,\tilde\rho}$ as independent indeterminates, we can 
now compute the permutation characters of $G(q)$ on $P_J(q)$ for various $J$.

\begin{Rem}\label{rempermprop}
We get constraints
on the possible values of the scalars $\zeta_{u,\tilde\rho}$ from general
facts about permutation characters.
More precisely:

\begin{itemize}
\item[(a)]
The case of the trivial character of $G(q)$ (the special case
$J=\Delta$, $P(q) = G(q)$) yields all $\zeta_{u,\tilde 1}$.
\item[(b)]
The values of permutation characters are non-negative rational
integers. 
From this we often  see that $\zeta_{u,\tilde\rho} \in \Q$, and so
$\zeta_{u,\tilde\rho} \in \{\pm 1\}$.
\item[(c)]
If only one of the remaining possibilities for one or several
$\zeta_{u,\tilde\rho}$ leads to non-negative
values then the scalars are determined.
\item[(d)]
The Green function $Q_w$ for the maximally split torus is the permutation
character ${1_{B(q)}}^{G(q)}$ restricted to unipotent elements. 
From its values we can find the sizes
of the intersections of the unipotent classes in $G(q)$ with $B(q)$ or
$U(q)$: For $u \in G(q)$ unipotent we have
\[ 
\begin{array}{rcl}
{1_{B(q)}}^{G(q)}(u) & = &\frac{1}{|B(q)|}\left| 
\{g \in G(q)\mid \; u^g \in U(q)\} \right| \\[3mm]
 & = &
| u^{G(q)} \cap B(q)| \; |C_{G(q)}(u)| \,/\, |B(q)|.
\end{array}
\]
All these values must be positive integers.
\end{itemize}
We can also use the general  facts about Green functions given 
in~\cite[7.6]{Ca85}:
\begin{itemize}
\item[(e)]
The values of the Green function $Q_{w_i}$ are rational integers.
\item[(f)] We have
$ \sum\limits_{u \in G(q) \textrm{ {\scriptsize unipotent}}} Q_{w_i}(u) =
\frac{|G(q)|}{|T_{w_i}(q)|}.$
\end{itemize}
We could also mention the orthogonality relations for Green functions, but
in our applications they never provided useful additional constraints.
\end{Rem}

\subsection{Permutation characters by counting fixed points}\label{ssecPcount}

To find additional equations for the scalars $\zeta_{u,\tilde\rho}$ we
compute some values of permutation characters by counting fixed points.
Note that for any finite group $G$ and subgroup $H \leq G$ the value
of the permutation character on the right cosets of $H$ for $g \in G$ is 
equal to
\[ {1_H}^G(g) = |\{Hx \mid\; x \in G, Hx = Hxg\}|. \]

To apply this to parabolic subgroups we need a set of representatives of
the right cosets of $P_J(q) \leq G(q)$. 

\begin{Pro}\label{repsmodP}
Let $J \subseteq \Delta$ be a subset of the simple roots of $\bG$ and 
$\bP_J$ be the corresponding standard parabolic subgroup and $W_J$ the
parabolic subgroup of the Weyl group generated by $J$. We say that
$w \in W$ is $J$-reduced  if $l(s_r w) > l(w)$ for all $r \in J$ (that is
$w$ is the shortest possible representative of the coset $W_J w \subseteq
W$).
\begin{itemize}
\item[(a)] A set of right coset representatives of $\bP_J\setminus \bG$ is
given by
\[ \{ \dot w \prod_{r \in \Phi^+,\, w(r) \in \Phi^-} x_r(a_r) \mid\;
w \textrm{ is $J$-reduced,}\; a_r \in \mathbf{k}\}.\]
\item[(b)] When $F$ acts trivially on $\Phi$ and $W$, then a set of right
coset representatives of $P_J(q)\setminus G(q)$ is given by
\[ \{ \dot w \prod_{r \in \Phi^+,\, w(r) \in \Phi^-} x_r(a_r) \mid\;
w \textrm{ is $J$-reduced,}\; a_r \in \Fq\}.\]
\end{itemize}
The factors in the products are taken in any fixed order.
\end{Pro}

Part~(b) can be generalized to twisted groups, but the statement becomes
more technical. We do not need this in the remainder of this article and
omit it.

\begin{Prf}
Note that for $J = \{\}$, $\bP_J = \bB$ the Borel subgroup, this is just a
reinterpretation of the strong form of the Bruhat decomposition. We consider
general $J$.

Recall that $\bP_J$ is generated by the $\bB \dot s_r \bB$ with $r \in J$.
Let $g' \in \bG$ and $w' \in W$ such that $g' \in \bB \dot w' \bB$. When $w'$ is
not $J$-reduced, there is $r \in J$ with $l(s_r w') < l(w')$ ($\Leftrightarrow
w'^{-1}(r) \in \Phi^-$ by~\cite[2.2.1]{Ca72}). We set $w = s_r w'$ and have
$(\bB \dot s_r \bB) (\bB \dot w \bB) = \bB \dot w' \bB$. This shows that the
coset $(\bP_J) g' = (\bP_J) g$ for some $g \in \bB \dot w \bB$. 
Since the length of
$w$  is smaller than the length of $w'$ we will find after a finite
number of applications of this step $g$ and $w$ such that $w$ is
$J$-reduced. Since $\bB \leq \bP_J$ we see from the strong form of the Bruhat
decomposition that the elements given in~(a) contain representatives of all
right cosets in $\bP_J\setminus \bG$. 

The multiplication rule for $(\bB \dot s_r \bB)(\bB \dot w \bB)$ shows that
for $g \in (\bB \dot w \bB)$, $g' \in (\bB \dot w' \bB)$ with
$(\bP_J)g  = (\bP_J) g'$ we have $(W_J) w = (W_J) w'$. Furthermore, 
for any $w' \in W$ there is exactly one $J$-reduced $w \in W$ with
$(W_J) w = (W_J) w'$ (the $J$-reduced $w$ maps $\Phi_J^+$ to itself
and only the trivial element of $W_J$ has this property).
This shows that  any $g' \in \bG$ determines a unique $J$-reduced $w$ such
that $(\bP_J)g  = (\bP_J) g'$  for some $g$ of the form $g = \dot w 
\prod\limits_{r \in \Phi^+,\, w(r) \in \Phi^-} x_r(a_r)$.

Finally, let $w \in W$ be $J$-reduced, $u, u' \in U_w$, and 
$(\bP_J)\dot w u = (\bP_J)\dot w u'$. The last condition is equivalent
to $\dot w u u'^{-1} \dot w^{-1} \in \bP_J$. We have $u u'^{-1} \in U_w$
and $\dot w U_w \dot w^{-1} = \prod\limits_{r \in \Phi^+,\, w(r) \in \Phi^-}
X_{w(r)}$. Since $w$ is $J$-reduced, so $w(\Phi_J^+) = \Phi_J^+$, we see
$\dot w U_w \dot w^{-1} \cap \bP_J = \{1\}$. Hence $u u'^{-1} = 1$ and $u =
u'$ and we conclude that the set of elements in~(a) contains a unique
representative for each right coset of $\bP_J$ in $\bG$.

The proof of~(b) is the same using the corresponding finite groups instead
of $\bG$, $\bB$, $\bP_J$ and $X_r$.
\end{Prf}
\medskip

\subsection{Computing a character value}\label{compval}
Let 
\[ v = x_{j_1}(c_1) \cdots x_{j_k}(c_k) \in U(q),\]
with $c_1, \ldots, c_k \in \Fq\setminus \{0\}$ be a 
unipotent element for which we want to
find the value of the permutation character ${1_{P_J(q)}}^{G(q)} (v)$.

For each $J$-reduced $w \in W$ we fix an ordering $r_1, \ldots, r_N$ of
the positive roots such that for some $l$ we have $w(r_i) \in \Phi^+$ for
$1 \leq i \leq l$ and $w(r_i) \in \Phi^-$ for $l+1 \leq i \leq N$.
We consider all cosets represented by
\[ \dot w x_{l+1}(a_{l+1}) \cdots x_N(a_N),\textrm{ with
} a_i \in \Fq \]
at once. For the computation we use independent indeterminates $y_{l+1}, 
\ldots, y_N$ over $\Fq$
(or even over $\Z$) instead of the $a_i$ in the expression above. 
We have to count for how many specializations 
$a_i \in \Fq$ of the $y_i$ the elements
\[ g_1 = \dot w x_{l+1}(a_{l+1}) \cdots x_N(a_N) \textrm{ and } g_2 = g_1 v\]
are in the same right coset of $P_J(q)$.
We apply the reordering algorithm in Proposition~\ref{reorder}(b) to 
rewrite
\[  \dot w \; x_{l+1}(y_{l+1}) \cdots x_N(y_N) x_{j_1}(c_1) \cdots 
x_{j_k}(c_k)\]
in the form
\[ \dot w \; x_1(b_1) \cdots x_{l}(b_l) x_{l+1}(b_{l+1}) \cdots x_N(b_N),\]
where we get $b_1, \ldots, b_N$ as polynomials in the indeterminates
$y_{l+1}, \ldots, y_N$.

The chosen ordering of the positive roots yields that any specialization of
the indeterminates in 
\[ \dot w \; x_1(b_1) \cdots x_{l}(b_l) \dot w^{-1} \]
in $\Fq$ yields an element in $U(q) \subset P_J(q)$. 

So, we want to count the tuples $(a_{l+1},\ldots, a_N) \in \Fq^{N-l}$ such
that specializing the $y_i$ to  $a_i$ in the expressions
\[ \dot w\; x_{l+1}(y_{l+1}) \cdots x_N(y_N) \textrm{ and }
\dot w\; x_{l+1}(b_{l+1}) \cdots x_N(b_N) \]
yields representatives of the same right coset of $P_J(q)$. 

Using Proposition~\ref{repsmodP}(b) this translates to counting the
solutions in $\F_q^{N-l}$ of the system of polynomial equations
\[ y_{l+1} - b_{l+1} = 0, \ldots, y_N - b_N = 0.\]

\subsection{Counting solutions}

We have reduced the computation of values of permutation characters to
counting solutions over $\Fq$ of systems of multivariate polynomial
equations over $\Fq$. In general, such counting is a difficult task. In
principle one could use Gröbner bases with respect to some elimination
monomial ordering. But computing Gröbner bases is itself a difficult
problem. We mention here some heuristics which often solve the problem for
systems of equations that occur in the context of the previous subsection.

\begin{itemize}
\item
For small $q$ we can reduce the degrees of the polynomials by the
substitution $y = y^q$ for any indeterminate $y$ (since $a^q-a = 0$ for all
$a \in \Fq$).
\item
It happens quite often that we have an equation of the form
$y = f$ with  $f$ a polynomial in indeterminates
different from $y$. We use this to substitute $y$ accordingly in all
equations, and hence eliminate $y$ from the system. We prefer cases where
$f$ has few terms, and we fix an upper bound for the amount of memory to use
and stop substitutions when the limit is reached.
\item 
There are many examples where the previous steps lead to an equation
$c = 0$ for a non-zero constant $c$ (so, there are no solutions) or that we
end up with no equation (so, the number of solutions is $q^k$ when there are
$k$ independent variables left).
\item
Otherwise, we use for fixed $q$ a straight forward backtrack search:
Specialize one
variable  to all possible values and solve for each value recursively 
the resulting system of equations in the remaining variables.
\end{itemize}

\subsection{Detecting cases with no solutions}\label{nocomp}
In the setup of Section~\ref{compval} we can often detect quickly that the
system of polynomial equations we get for a $J$-reduced $w \in W$ has no
solutions.

We consider the roots $\Psi' = \{r_{j_1}, \ldots, r_{j_k}\}$ 
corresponding to the
non-trivial root element factors of the element $v$ for which we want to
compute the value of a permutation character. We 
determine the subset $\Psi \subseteq \Psi'$ consisting of $r$ 
which correspond to only one factor in the given product $v$ and such
that $r' \not\prec r$ for all $r' \in \Psi' \setminus \{r\}$.

Now let $w \in W$ be  $J$-reduced and assume that $w(r) \in \Phi^-$ for some
$r \in \Psi$. Let $r_i =r $ for some $l+1 \leq i \leq N$ and $x_{r_i}(c)$
(with $c \neq 0$) be
the corresponding factor in $v$. Then we see from Proposition~\ref{reorder}(c) 
that  our method leads to  an equation $y_i = y_i + c$ causing that there are
no solutions.

So, we can skip our computation for the double coset of $w$ whenever there
is at least one $r \in \Psi$ such  that $w(r)$ is negative.

\section{Application to $\mathbf{E_8(q)}$}\label{secE8}

As an application we want to determine the (ordinary) Green functions for
the exceptional groups $\bG^F = G(q) = E_8(q)$ in 
bad characteristic $p = 2, 3, 5$.
For good characteristic $p > 5$ this problem was solved by Beynon and
Spaltenstein in~\cite{BS84}. In Section~\ref{secGFalg} we have explained
that we need to find appropriate class representatives of the unipotent
classes in $\bG$ and certain associated complex scalars $\zeta_{u,\tilde\rho}$
of absolute value $1$. The following theorem summarizes the result.

\begin{Thm}\label{thm}
With one exception we can find in each unipotent conjugacy class $C$ of
$E_8(\barGF{q})$ an element $u \in C^F$ such that $F$ acts trivially on
$A(u)$ and such that we have for all associated
scalars $\zeta_{u,\rho} = 1$.

The exception is the class $C$ labeled $D_8(a_3)$ 
when $q \equiv -1 (3)$ (so $p \neq 3$). In this case there is a unique
$\bG^F$-class of $u \in C^F$ such that $F$ acts trivially on $A(u) \cong S_3$.
Then $\zeta_{u,\rho} = 1$ except when $\rho = -1$ is the sign character 
of $S_3$ where $\zeta_{u,-1} = -1$.
\end{Thm}

The starting point of the proof is Spaltenstein's table in~\cite{Sp85} which
describes the generalized Springer correspondence for groups of type $E_8$
(the description is slightly different in the cases $p=2,3,5,>5$). The table
has labels for the unipotent classes in the algebraic group and also gives
the dimensions of the classes and the isomorphism types of the component
groups $A(u)$. 

The unipotent conjugacy classes were determined by Mizuno~\cite{Mi80} and
Spaltenstein used the labeling of Mizuno. We will use some explicit class
representatives found by Mizuno. Some care is needed when using Mizuno's
tables because they contain (very few) errors. An important correction of
two $A(u)$ in characteristic $2$ is mentioned in~\cite[5.5]{Sp85}.
For some explicit computations in the groups we want to use elements given
in Mizuno's paper. Therefore, we use the same structure constants as Mizuno
and read off the extraspecial signs he used (there are some errors in the
table of structure constants in the paper, but we recompute everything from
the signs for extraspecial pairs.)

We can conclude from Mizuno's results that all unipotent classes contain an
$F$-stable element $u$ such that $F$ acts trivially on $A(u)$: Otherwise
there would be a class $C$ such that $F$ acts on $A(u)$ as a non-inner
automorphism for any $u \in C^F$. All groups occurring as $A(u)$ have the
easy to check property that the number of $\phi$-conjugacy classes for any
non-inner automorphism $\phi$ is smaller than the number of conjugacy classes.
So, for some power $F^m$ of $F$ such that $F^m$ acts trivially on all $A(u)$
the group $\bG^{F^m}$ would have more unipotent classes than $\bG^F$. 
But Mizuno showed
that the number of unipotent classes only depends on the characteristic $p$.

From now on we assume that $p = 2$. The cases $p=3,5$ must be considered
separately, but the arguments to obtain our result are the same (and
sometimes a bit easier because there are fewer unipotent classes).

\textbf{Step (1).}
So, assuming that we choose representatives $u \in C^F$ with trivial action
of $F$ on $A(u)$ for each unipotent class $C$ of $\bG$, we can use the
Lusztig-Shoji algorithm to compute the (generalized) Green functions as
linear combinations of the functions $\zeta_{u,\rho} Y_{u,\rho}$, see
Section~\ref{secGFalg}. For this we only need the information in
Spaltenstein's table and the character tables of (relative) Weyl groups,
which are, e.g., available in GAP~\cite{GAP} or CHEVIE~\cite{CHEVIE}.
For the computations in this step we use independent indeterminates instead
of the so far unknown scalars $\zeta_{u,\rho}$.

Up to this point we do not need to be more specific on the choice of class
representatives $u \in C^F$. (But, of course, the $\zeta_{u,\rho}$ depend
on the choice of $u$.)

We have 74 unipotent classes in $\bG$ and 146 unipotent classes in $G(q) =
E_8(q)$ (where $q$ is any power of $p=2$). We need to determine 112 scalars
$\zeta_{u,\rho}$ which so far are represented by 112 indeterminates in our
table of Green functions.

\textbf{Step (2).}
We consider the trivial character of $G(q)$ (the permutation character on
$P_J(q)$ for $J=\Delta$) and use Remark~\ref{rempermprop}(a). This shows
that the 74 scalars $\zeta_{u,1}$ all must be $1$ (independent of the choice
of representatives $u$). We specialize their corresponding indeterminates in
our table of Green functions such that only 38 unknown $\zeta_{u,\rho}$ remain.

From the rationality of Green functions, see Remark~\ref{rempermprop}(e)
or~(b), we see that all of these $\zeta_{u,\rho} \in \{\pm 1\}$.

\textbf{Step (3).}
We determine the sizes of the unipotent classes in $G(q)$ with the
information from Remark~\ref{remLambda}. In many cases they are 
uniquely determined from the matrix $\Lambda$, and they are the same for all
possible values of the not yet known $\zeta_{u,\rho}$.
For example this is the case for the class $D_8(a_3)$ (where $A(u) \cong
S_3$ for $u$ in this class).

In other cases certain possibilities for $\zeta_{u,\rho}$ do not lead to 
a diagonal matrix of class lengths as described in~\ref{remLambda}. For
example for the class $D_4(a_1)$, also with $A(u) \cong S_3$, we get
that the two not yet known scalars can only be $1$, and this determines the
class lengths.

An interesting case are classes $C$ where for $u \in C$ we have $A(u)$ of
order $2$ and the two classes in $C^F$ have different lengths. There both 
choices in $\{\pm 1\}$ of the scalar for the non-trivial character 
of $A(u)$ lead to a diagonal matrix with both class lengths, but in different
order. So, here we can just assume that the scalar is $\zeta_{u,-1} = 1$
and we get the class length of $u$.

Another case is the class $D_6(a_1)$ with $A(u)$ elementary abelian of
order $4$, where only one unknown scalar appears in the ordinary Green
functions. We get the class lengths but not the scalar. Also for the class
$2A_4$ with $A(u)\cong S_5$ we get the class lengths but not the scalars.

Now we use the property of Green functions from Remark~\ref{rempermprop}(f)
for the split torus of order $(q-1)^8$. We can evaluate this sum using the
lengths of the unipotent classes. 

We get an equation of the form
\[ z_1 f_1(q) + \ldots + z_{18} f_{18}(q) = \frac{|G(q)|}{(q-1)^8}. \]
where $z_1, \ldots, z_{18}$ are 18 out of our 38 indeterminates remaining
after step~(2)
and the $f_1(q), \ldots, f_{18}(q)$ are polynomial expressions in $q$.

For all $f_i(q)$ it is easy to see that they evaluate to positive integers
for $q$ any power of $2$. Also, the sum of the $f_i(q)$ equals the right
hand side of the equation. So, the only way to satisfy the equation by
substituting all $z_i$ by complex numbers of absolute value $1$ is to set
all $z_i$ to $1$.

\textbf{Step (4).}
Let $C$ be one of the classes $(A_5+A_1)''$, $A_5+2A_1$, $D_6+A_1$, $D_8$ 
where for $u\in C$ we have $|A(u)| = 2$. In these cases the two
classes in $C^F$ have the same order and the two possibilities
for $\zeta_{u,-1} \in \{\pm 1\}$ exchange all values of Green functions on
these two classes. So, there is a choice of $u \in C^F$ such that
$\zeta_{u,-1} = 1$.

A similar argument works for the class $C = E_7+A_1$ where for $u \in C$ the
group $A(u)$ is elementary abelian of order $4$. Setting the unknown scalar
to $\{\pm 1\}$ leads to the same permutation of the values 
on the four classes in $C^F$ for all Green functions. 
So, there is a choice of $u$ such
that the scalar is $1$.

We have 15 remaining indeterminates.

\textbf{Step (5).}
Now we use Remark~\ref{rempermprop}(d) and compute the intersections of
unipotent classes with the Borel subgroup $B(q)$. Some of these expressions
are rational polynomials in $q$ and one or several of the unknown scalars.

The expressions are big, but now we also use Proposition~\ref{proprimefield}
and specialize $q$ to $q=2$ to get some rational linear combinations of
unknown scalars which must be positive integers. For example, the
expressions for the class $A_5+A_2$ with component group $S_3$ contain two
unknown scalars. Setting both scalars to $1$ yields positive integers,
but the other three possibilities either yield a negative or a non-integer
rational number. (Note that in the case of the non-abelian $A(u) \cong S_3$ 
the $G(q)$-class of $u$ is already fixed by the condition that $F$ acts 
trivially on $A(u)$.)

With this method we find another 11 of the scalars to be $1$. 

\textbf{Step (6).}
Two of the remaining four unknown scalars belong to the class $C=D_8(a_3)$.
For $u \in C$ we have $A(u) \cong S_3$, so that the $G(q)$-class 
of $u$ is fixed by the condition that $F$ acts trivially on $A(u)$.
We still need to determine $\zeta_{u,-1}$ and $\zeta_{u,\rho}$ where we
write $-1$ for the sign character and $\rho$ for the irreducible 
character of degree $2$ of $S_3$.

Computing values of permutation characters on unipotent elements for various
parabolic subgroups $P_J(q)$ as in Section~\ref{ssecPDL} 
we find that, e.g.,  for $J$ of type $D_4$ the mentioned scalars appear in the 
values on the classes in $C^F$. 

Using Proposition~\ref{proprimefield} we specialize $q = 2$. Then the four
possibilities for the scalars would lead to different character values
of ${1_{P_J(q)}}^{G(q)}$ on the classes in $C^F$, 
namely the tuples: $(92897, 177889, 89825 )$,
$(179937, 90849, 176865)$, $(88801, 177889, 91873)$, $(175841, 90849,
178913)$.

Now we use Section~\ref{compval} and 
compute the values of this permutation character on the elements
$z_{77}, z_{78}, z_{79}$ given in Mizuno~\cite{Mi80} and find the first of
the tuples above which corresponds to $\zeta_{u,-1} = -1$ and
$\zeta_{u,\rho} = 1$. This is the exception stated in our Theorem~\ref{thm}.
Note that by Proposition~\ref{proprimefield} this exception only occurs when
$q$ is an odd power of $2$ (or $q \equiv -1 (3)$).

Remarks: The computation of each value took about 5 seconds, there are
$3628800$ $J$-reduced elements $w$, but only about $500$ of them do not
fulfill the criterion in Section~\ref{nocomp}.
Note that with our computation we can actually show that the elements 
$z_{77}, z_{78}, z_{79}$ in Mizuno's list are indeed in the class $C=D_8(a_3)$
because the same values do not occur on other classes (and we see that $F$ acts
trivially on the component group of $z_{77}$). In fact, we found
an error in Mizuno's list by computing
the value of $z_{64}$. This also yields the value $177889$ and
this shows that this element is conjugate to $z_{78}$.

\textbf{Step (7).}
The last two unknown scalars correspond to the class $C=D_8(a_1)$. For
$u \in C$ the group $A(u)$ is a dihedral group of order $8$ which has $5$ 
conjugacy classes and two non-trivial linear
characters $\epsilon', \epsilon''$ which appear in the Springer
correspondence. All values of (ordinary) Green functions are the same on the
classes of $u = u_1$ and $u_a$ where $a$ is the non-trivial element in the
center of $A(u)$.

As in step~(6) we can use the parabolic subgroup for $J$ of
type $D_4$. The four possibilities for the two unknown scalars always lead
to the same four character values, but they differ in which of the four
values appears twice. For $q=2$ the four values are $\{6785, 2625, 6401,
2241\}$. Using the five representatives $z_i$, $i \in \{44,46,47,48,49\}$, in
Mizuno's list and computing the character values for these elements we
get the value $6785$ twice, namely for $z_{44}$ and $z_{49}$. This shows
that also the scalars $\zeta_{u,\epsilon'} = \zeta_{u,\epsilon''} = 1$.

\textbf{Remarks.} In our first proof of the Theorem~\ref{thm}
we used arguments as in
steps~(6) and~(7), that is computations as explained in
Section~\ref{ssecPcount}, for many more of the unknown scalars. Even large
$J$ for which the computations are pretty fast provide a lot of useful
information. We were a bit surprised that in the revised version presented
here we could  avoid almost all of these computations. Nevertheless, as
hinted in the remarks to step~(6), the computational
method can be useful to obtain more detailed information as in
Theorem~\ref{thm}. For example, we have seen above that we can identify 
the $G(q)$-class of some unipotent elements. The values of some permutation
characters are useful class invariants which often distinguish classes from
all others.

Of course, using the techiques of this paper it is possible to check and
recover known results for other types of groups of small rank (where the
original argument involved much more elaborate computations).


\newcommand{\etalchar}[1]{$^{#1}$}

\end{document}